\documentclass[12pt,twoside,english,a4paper]{article}
\usepackage{babel,amssymb,amsthm,amsmath}
\usepackage{graphicx}
\usepackage{mathabx}
\usepackage{bm}

\newcommand{\punto}{\,\cdot\,}
\newcommand{\ds}{\displaystyle}
\newcommand{\smallfrac}[2]{{\textstyle\frac{#1}{#2}}}

\newcommand{\Nelt}{{N_{\mathrm{elt}}}}
\newcommand{\Nver}{{N_{\mathrm{ver}}}}
\newcommand{\Nnd}{{N_{\mathrm{qd}}}}
\newcommand{\Nqd}{{N_{\mathrm{qd2}}}}

\newcommand{\Nfc}{{N_{\mathrm{fc}}}}
\newcommand{\Ndir}{{N_{\mathrm{dir}}}}
\newcommand{\Nneu}{{N_{\mathrm{neu}}}}

\numberwithin{equation}{section}

\title{Matlab tools for HDG in three dimensions}

\date{\today}

\author{Zhixing Fu\footnote{Department of Mathematical Sciences, University of Delaware, USA. {\tt zfu@math.udel.edu}}, Luis F. Gatica\footnote{Departamento de Matem\'atica y F\'\i sica Aplicadas, Universidad Cat\'olica de la Sant\'\i sima Concepci\'on \& CI$^2$MA--Universidad de Concepci\'on, Concepci\'on, Chile. {\tt lgatica@ucsc.cl}. Research partially funded by Direcci\'on de Investigaci\'on, UCSC.}, Francisco--Javier Sayas\footnote{Department of Mathematical Sciences, University of Delaware, USA. {\tt fjsayas@math.udel.edu}. Research partially funded by the National Science Foundation (NSF-DMS 1216356)}}

\begin{document}

\maketitle

\begin{abstract}
In this paper we provide some Matlab tools for efficient vectorized coding of the Hybridizable Discontinuous Galerkin for linear variable coefficient reaction-diffusion problems in polyhedral domains. The resulting tools are modular and include enhanced structures to deal with convection-diffusion problems, plus several projections and a superconvergent postprocess of the solution. Loops over the elements are exclusively local and, as such, have been parallelized.
\end{abstract}

\section{Introduction}

In this paper we provide some programming tools for full Matlab implementation of the Hybridizable Discontinuous Galerkin (HDG) method on general conforming tetrahedral meshes for fixed but arbitrary polynomial degree. The presentation is detailed on a second order linear reaction-diffusion equation with variable coefficients and mixed boundary conditions, but we also provide the tools to construct the matrices needed for convection-diffusion problems with variable convection, thus creating all necessary blocks to deal with general steady-state problems.

The HDG method originated in a sequence of papers of Bernardo Cockburn and his collaborators, consolidating in the unified framework of \cite{CoGoLa:2009}. As seen in that paper, the HDG can be considered as a Mixed Finite Element Method (MixedFEM) \cite{BrFo:1991}, coded with the use of Lagrange multipliers to  weakly enforce the restrictions on interelement faces \cite{ArBr:1985}, and then hybridized so that the only global variable is the collection of Lagrange multipliers, that ends up being an optimal approximation of the primal variable on the faces of the triangulation. 

Compared with general MixedFEM (programmed in hybridized form), HDG has the advantage of not using degrees of freedom to stabilize the discrete equations, while keeping equal optimal order of convergence in all computed fields. (Stability is obtained through a stabilization parameter.) From this point of view, HDG is a valid option if one is willing to pay the prize of using MixedFEM, for instance, to obtain approximations of more fields of the solution of the problem. In comparison with other FEM, that work directly on the second order formulation, HDG performs well for high orders \cite{KiSpCo:2012}. For low order methods, HDG can be adopted in situations where either MixedFEM or DG are thought to be advantageous. We will, however, not exploit here the advantages of HDG/DG for having non-conforming meshes or variable degree. Extension of the code to variable-degree methods does not seem to change much, but it requires rethinking the data structures and the vectorization process. Extension to more general meshes would require new tools for the geometric handling that we are not dealing with at this moment.
In comparison to other DG methods (mainly those of the Interior Penalty family), HDG requires less degrees of freedom in the solution of the global system, since this has been reduced to the interfaces of the elements. On the positive side as well, HDG does not contain any penalization parameter that needs tuning to obtain convergence. It also has some attractive superconvergence properties that allow for local element-by-element postprocessing {\em \`a la Stenberg} \cite{Stenberg:1991}.

One goal of the paper is the systematization of the construction of local and global matrices by looping over quadrature nodes and polynomial degrees, avoiding large loops over elements. We partially accomplish this by using Matlab inbuilt functions for Kronecker products, construction of sparse matrices, and vectorization. All loops on elements are purely local and have been parallelized so that they can take advantage of the Matlab Parallel Toolbox. We hope this piece of work will contribute to the popularization of a method that has already a sizeable follow-up, given its good properties. This being Matlab code, we are not expecting the code to run on very large problems, but, as we show in the experiments, we can show reasonably high order of convergence  for three dimensional problems, working on a laptop with only two processors. We will also comment on how this code can be easily modified to provide the hybridized implementation of the Brezzi--DOuglas--Marini (BDM) mixed element, which is an alternative to using $\mathbf H(\mathrm{div})$-conforming bases \cite{Ervin:2012}.

\paragraph{Model equations.}
Let $\Omega\subset \mathbb R^3$ be a polyhedron with boundary $\Gamma$, divided into a Dirichlet and a Neumann part ($\Gamma_D$ and $\Gamma_N$) such that each part is the union of faces of $\Gamma$. The unit outward-poiting normal vector field on $\Gamma$ is denoted $\boldsymbol\nu$.
The model problem we will be discussing in this document is
\begin{subequations}\label{eq:1.1}
\begin{alignat}{6}
\label{eq:1.1a}
\kappa^{-1} \boldsymbol q+\nabla u &=0 &\qquad & \mbox{in $\Omega$},\\
\label{eq:1.1b}
\nabla\cdot\boldsymbol q+c\,u &=f & & \mbox{in $\Omega$},\\
\label{eq:1.1c}
u & = u_D & & \mbox{on $\Gamma_D$},\\
\label{eq:1.1d}
-\boldsymbol q\cdot\boldsymbol\nu &=\boldsymbol g_N\cdot\boldsymbol\nu & & \mbox{on $\Gamma_N$}.
\end{alignat}
\end{subequations}
The diffusion coefficient $\kappa$ is strictly positive, while $c\ge 0$. (Both of them are functions of the space variables.) The Neumann boundary condition is given in non-standard way, as the normal component of a vector field (of which only the normal component is used), in order to give an easier way to test exact solutions. Modifications for the case of a scalar field are straightforward.

\paragraph{Discrete elements.} The Hybridizable Discontinuous Galerkin method that we will describe here is based on regular tetrahedrizations of the domain. We thus consider a tetrahedral partition of $\Omega$, $\mathcal T_h$. The set of all faces in the triangulation is denoted $\mathcal E_h$, with the subsets $\mathcal E_h^{\mathrm{int}}$, $\mathcal E_h^D$, $\mathcal E_h^N$ corresponding to interior, Dirichlet and Neumann faces. For convenience, we will write
\[
\Nelt=\# \mathcal T_h, \qquad \Nfc=\# \mathcal E_h, \qquad \Ndir=\# \mathcal E_h^D, \qquad \Nneu=\#\mathcal E_h^N.
\]
Upon numbering, we can identify elements and faces with respective index sets
\[
\mathcal T_h\equiv \{1,\ldots,\Nelt\}, \qquad \mathcal E_h\equiv \{1,\ldots,\Nfc\}.
\]
This will allow us to write some computational expressions in a format that is very close to their mathematical definition.
The local spaces for discretization of $u$ and $\boldsymbol q$ are those of trivariate polynomials of degree up to $k$. The global description of these spaces is
\[
 W_h:=\prod_{K\in \mathcal T_h}\mathcal P_k(K), \qquad \boldsymbol V_h:=W_h^3=\prod_{K\in \mathcal T_h} \mathcal P_k(K)^3.
\]
There is a third space, defined on the skeleton of the tetrahedrization:
\[
M_h:=\prod_{e\in \mathcal E_h} \mathcal P_k(e),
\]
where $\mathcal P_k(e)$ is the space of bivariate polynomials of degree not larger than $k$ on tangential coordinates. Integral notation will always be given as
\[
(u,v)_K:=\int_K u\,v,\qquad (\boldsymbol q,\boldsymbol r)_K:=\int_K \boldsymbol q\cdot\boldsymbol r, \qquad \langle u,v\rangle_{\partial K}:=\int_{\partial K}u\, v.
\]
Similarly, we will use terms like $\langle u,v\rangle_e$ with $e\in \mathcal E_h$, $\langle u,v\rangle_{\Gamma_D}$, and $\langle u,v\rangle_{\Gamma_N}$. On the boundary of a given element, the normal vector $\boldsymbol\nu_K$ will point outwards. However, when it is clear what the element is, we will simply write $\boldsymbol\nu$.

\paragraph{HDG.} A key ingredient of HDG is a stabilization function. This is a non-negative piecewise constant function on the boundary of each triangle. Thus,
\[
\tau_{K}|_e \in \mathcal P_0(e) \quad \forall e\in \mathcal E(K), \quad \forall K\in \mathcal T_h, \qquad\mbox{and}\qquad \tau_{K}\ge 0.
\]
Here $\mathcal E(K)=\{ e^K_1,\ldots,e^K_4\}$ is the ordered set of faces of $K$. The function $\tau$ is not single-valued in internal faces. We demand that for each $K$, the function $\tau_{K}$ cannot vanish identically, that is, there exists $e\in \mathcal E(K)$ such that $\tau_{K}|_e>0$. 
The HDG method works separately each of the equations in \eqref{eq:1.1}. There are three unknowns: $(\boldsymbol q_h,u_h,\widehat u_h)\in \boldsymbol V_h\times W_h\times M_h$. Locally, we will think of $(\boldsymbol q_K,u_K)\in \mathcal P_k(K)^3\times \mathcal P_k(K)$, while $\widehat u_h$ will be counted by global face numbering, so we will sometimes refer to values $\widehat u_e$. 
\begin{subequations}\label{eq:1.2}
The PDE \eqref{eq:1.1a}-\eqref{eq:1.1b} is discretized element by element with the equations
\begin{alignat}{4}
\label{eq:1.2a}
(\kappa^{-1}\boldsymbol q_K,\boldsymbol r)_K - (u_K,\nabla\cdot\boldsymbol r)_K+\langle \widehat u_h,\boldsymbol r\cdot\boldsymbol\nu_K\rangle_{\partial K} &=0 &\qquad &\forall\boldsymbol r\in \mathcal P_k(K)^3,  \\
\label{eq:1.2b}
(\nabla\cdot\boldsymbol q_K,w)_K+(c\,u_K,w)_K+\langle \tau_K(u_K-\widehat u_h),w\rangle_{\partial K} &=(f,w)_K & &\forall w\in \mathcal P_k(K),
\end{alignat}
for all $K\in \mathcal T_h$. All the remaining equations are defined (and counted) on edges. We first have flux equilibrium on internal faces: for all $\mathcal E_h^{\mathrm{int}}\ni e =K\cap \widetilde K$, we impose
\begin{equation}\label{eq:1.2c}
\langle \boldsymbol q_K\cdot\boldsymbol\nu_K+\tau_K(u_K-\widehat u_e),\widehat v\rangle_e+\langle \boldsymbol q_{\widetilde K}\cdot\boldsymbol\nu_{\widetilde K}+\tau_{\widetilde K}(u_{\widetilde K}-\widehat u_e),\widehat v\rangle_e=0 \quad \forall \widehat v\in \mathcal P_k(e).
\end{equation}
This means that the normal numerical flux, $\Phi_K:= -\boldsymbol q_K\cdot\boldsymbol\nu_K-\tau_K(u_K-\widehat u_e)$, which is an element of $\mathcal P_k(e)$ for all $e\in \mathcal E(K)$, is essentially single valued on internal faces, that is $\Phi_K+\Phi_{\widetilde K}=0$ on $e=K\cap \widetilde K$. We finally impose the boundary conditions on Dirichlet faces $e\in \mathcal E_h^D$
\begin{equation}\label{eq:1.2d}
\langle \widehat u_e,\widehat v\rangle_e=\langle u_D,\widehat v\rangle_e\qquad \forall \widehat v\in \mathcal P_k(e),
\end{equation}
and on Neumann faces $e\in \mathcal E_h^N$, 
\begin{equation}\label{eq:1.2e}
-\langle \boldsymbol q_K\cdot\boldsymbol\nu_K+\tau_K(u_K-\widehat u_e),\widehat v\rangle_e=\langle\boldsymbol g_N\cdot\boldsymbol\nu_K,\widehat v\rangle_e \qquad \forall \widehat v\in \mathcal P_k(e), \mbox{ where } e\in\mathcal E(K).
\end{equation}
\end{subequations}
Equations \eqref{eq:1.2} make up for a square system of linear equations. The hybridization of the methods is a static condensation (substructuring) strategy that allows to write the method as a system of equations where only $\widehat u_h$ appears as an unknown. For more methods that fit in this framework, see \cite{CoGoLa:2009}. Theory has been developed in a series of papers, but revisited and deeply reorganized in \cite{CoGoSa:2010}. Readers acquainted with programming mixed finite element methods will recognized the hybridized form of \cite{ArBr:1985}. (We will come back to this at the very end of the paper.)

\section{Geometric structures}\label{sec:2}

\paragraph{Tetrahedra.} All elements will be mapped from the reference tetrahedron $\widehat K$ with vertices
\[
\widehat{\mathbf v}_1:=(0,0,0) \qquad \widehat{\mathbf v}_2:=(1,0,0), \qquad \widehat{\mathbf v}_3:=(0,1,0), \qquad \widehat{\mathbf v}_4:=(0,0,1).
\]
Note that $|\widehat K|:=\mathrm{vol}\,\widehat K=1/6.$ Given a tetrahedron with vertices $(\mathbf v_1,\mathbf v_2,\mathbf v_3,\mathbf v_4)$ (the order is relevant), we consider the affine mapping $F_K:\widehat K\to K$ 
\begin{equation}\label{eq:2.0}
F_K(\widehat{\mathbf x})=\mathrm B_K \widehat{\mathbf x}+\mathbf v_1, \qquad 
\mathrm B_K=
\left[\begin{array}{ccc}
x_2-x_1 & x_3-x_1 & x_4-x_1\\ 
y_2-y_1 & y_3-y_1 & y_4-y_1\\ 
z_2-z_1 & z_3-z_1 & z_4-z_1
\end{array}\right].
\end{equation}
This map satisfies
$F_K(\widehat{\mathbf v}_i)=\mathbf v_i$, $ i\in\{1,2,3,4\}.$
All elements of the triangulation will be given with positive orientation, that is,
\[
\mathrm{det}\,\mathrm B_K=\Big( (\mathbf v_2-\mathbf v_1)\times (\mathbf v_3-\mathbf v_1)\Big)\cdot (\mathbf v_4-\mathbf v_1)>0.
\]
In this case $\mathrm{det}\,\mathrm B_K=6\,|K|$.

\paragraph{Faces.} A triangle $e$ in $\mathbb R^3$ with vertices $(\mathbf w_1, \mathbf w_2, \mathbf w_3)$ (the order is relevant), will be parametrized with 
\[
\boldsymbol\phi_e(s,t):=s\,(\mathbf w_2-\mathbf w_1)+t\,(\mathbf w_3-\mathbf w_1)+\mathbf w_1, \quad \boldsymbol\phi_e:\widehat K_2:=\{(s,t)\,:\,s,t\ge 0, s+t\le 1\}\to e.
\]
We note that $|\partial_s\boldsymbol\phi_e\times \partial_t\boldsymbol\phi_e|=2|e|$, where $|e|$ is the area of $e$. The local orientation of the vertices of $e$ gives an orientation to the normal vector. We will define the normal vector so that its norm is proportional to the area of $e$, that is
\[
\mathbf n_e:=\smallfrac12 \Big( (\mathbf w_2-\mathbf w_1)\times (\mathbf w_3-\mathbf w_1)\Big).
\]
Also, if  $\widehat{\mathbf w}_1:=(0,0)$, $\widehat{\mathbf w}_2:=(1,0)$, $\widehat{\mathbf w}_3:=(0,1),$
then $\boldsymbol\phi_e(\widehat{\mathbf w}_i)=\mathbf w_i$, for $i\in \{1,2,3\}$.

\paragraph{Boundaries of the tetrahedra.} Given a tetrahedron $K$ with vertices $(\mathbf v_1,\mathbf v_2,\mathbf v_3,\mathbf v_4)$ we will consider its four faces given in the following order (and with the inherited orientations):
\begin{equation}\label{eq:2.1}
\begin{array}{ccc}
e^K_1  &\qquad \longleftrightarrow \qquad& (\mathbf v_1,\mathbf v_2,\mathbf v_3)\\[1.3ex]
e^K_2  &\qquad \longleftrightarrow \qquad& (\mathbf v_1,\mathbf v_2,\mathbf v_4)\\[1.3ex]
e^K_3  &\qquad \longleftrightarrow \qquad& (\mathbf v_1,\mathbf v_3,\mathbf v_4)\\[1.3ex]
e^K_4  &\qquad \longleftrightarrow \qquad& (\mathbf v_4,\mathbf v_2,\mathbf v_3).
\end{array}\qquad\qquad \left[\begin{array}{ccc} 1 & 2 & 3 \\[1.3ex] 1 & 2 & 4 \\[1.3ex]1 & 3 & 4 \\[1.3ex]4 & 2 & 3\end{array}\right]
\end{equation}
(Note that with this orientation of the faces, the normals of the second and fourth faces point outwards, while those of the first and third faces point inwards. This numbering is done for the sake of parametrization.) For integration purposes on $\partial K$, we will use parametrizations of the  faces $ e^K_\ell \in \mathcal E(K)$
\[
\boldsymbol\phi_\ell^K:\widehat K_2\to e^K_\ell \qquad \ell\in \{1,2,3,4\},
\]
given by the formulas
\begin{equation}\label{eq:2.5}
\begin{array}{l}\ds \boldsymbol\phi_1^K(s,t) :=F_K (s,t,0),\\[1.5ex] 
\ds \boldsymbol\phi_2^K(s,t) :=F_K(s,0,t),\\[1.5ex] 
\ds \boldsymbol\phi_3^K(s,t) :=F_K(0,s,t), \\[1.5ex] 
\ds \boldsymbol\phi_4^K(s,t) :=F_K(s,t,1-s-t).\end{array}
\end{equation}
Consider the affine invertible maps $F_\mu:\widehat K_2\to\widehat K_2$ given by the formulas
\begin{equation}\label{eq:2.2}
\begin{array}{l}
\ds F_1(s,t) := (s,t)\\[1.5ex] 
\ds F_2(s,t) := (t,s)\\[1.5ex] 
\ds F_3(s,t) := (t,1-s-t)\\[1.5ex] 
\ds F_4(s,t) := (s,1-s-t)\\[1.5ex] 
\ds F_5(s,t) := (1-s-t,s)\\[1.5ex] 
\ds F_6(s,t) := (1-s-t,t)\\[1.5ex] 
\end{array} \qquad \qquad 
\left[\begin{array}{ccc} \mathbf 1 & \mathbf 2 & \mathbf 3 \\[1.5ex] \mathbf 1 & 3 & 2 \\[1.5ex] 3 & 1 & 2 \\[1.5ex] 3 & \mathbf 2 & 1 \\[1.5ex] 2 & 3 & 1 \\[1.5ex] 2 & 1 & \mathbf 3\end{array}\right].
\end{equation}
The table on the right shows the indices of the images of the vertices $(\widehat{\mathbf w}_1,\widehat{\mathbf w}_2,\widehat{\mathbf w}_3)$, with boldface font for those that stay fixed. We note that $F_2$, $F_4$ and $F_6$ change orientation. Take now a tetrahedron $K$, and assume that $e=e^K_\ell$, i.e., $e^K_\ell$  is the face $e$ in a global list of faces. We thus have six possible cases 
of how the parametrizations $\boldsymbol\phi_\ell^K$ and $\boldsymbol\phi_e$ match. We will encode this information in a matrix $\mathrm{perm}(K,\ell)$ so that
\begin{equation}\label{eq:2.3}
\boldsymbol\phi_e\circ F_\mu=\boldsymbol\phi_\ell^K, \mbox{ if } e=e^K_\ell \mbox{ and } \mu=\mathrm{perm}(K,\ell).
\end{equation}
We will refer to this matrix as the permutation matrix.

\paragraph{Data structure.} The basic tetrahedrization of $\Omega$ (including information on Dirichlet and Neumann boundaries) is given through four fields of a data structure {\tt T}:
\begin{itemize}
\item {\tt T.coordinates} is an $\Nver\times 3$ matrix with the coordinates of the vertices of the triangulation.
\item {\tt T.elements} is an $\Nelt\times 4$ matrix, whose $K$-th row of the matrix contains the indices of the vertices of $K$.
\item {\tt T.dirichlet} is an $\Ndir\times 3$ matrix, with the vertex numbers for the Dirichlet faces.
\item {\tt T.neumann} is an $\Nneu\times 3$ matrix, with the vertex numbers for the Neumann faces.
\end{itemize}
Positive orientation, of listings of vertices for elements and faces, is always assumed, 
In {\em expanded form}, the tetrahedral data structure contains many more useful fields. All these elements can be easily precomputed. It will be useful for what follows to assume that they are easy to access whenever needed.
\begin{itemize}
\item {\tt T.faces} is an $\Nfc\times 4$ matrix with a list of faces: the first three columns contain the global vertex numbers for the faces (its order will give the intrinsic parametrization of the face); Dirichlet and Neumann faces are numbered exacly as in {\tt T.dirichlet} and {\tt T.neumann}, the fourth column contains an index:
\begin{itemize}
\item 0 for interior faces
\item 1 for Dirichlet faces
\item 2 for Neumann faces
\end{itemize}
\item {\tt T.dirfaces} and {\tt T.neufaces} are row vectors with the list of Dirichlet and Neumann faces, that is, they point out what rows of {\tt T.faces} contain a $1$ (resp a $2$) in the last column.
\item {\tt T.facebyele} is an $\Nelt\times 4$ matrix, whose $K$-th row contains the numbers of faces that make up $\partial K$, with the faces given in the order shown in the table in \eqref{eq:2.1}. Note that this is the matrix we have described as $e^K_\ell$.
\item {\tt T.perm} is an $\Nelt\times 4$ matrix containing numbers from $1$ to $6$. Its $K$-th row indicates what permutations are needed for each of the faces to get to the proper numbering of the face, i.e., this is just the matrix $\mathrm{perm}(K,\ell)$.
\item {\tt T.volume} is an $\Nelt\times 1$ column vector with the volumes of the elements.
\item {\tt T.area} is an $\Nfc\times 1$ column vector with the areas of the faces.
\item {\tt T.normals} is an $\Nelt\times 12$ matrix with the {\em non-normalized} normal vectors for the faces of the elements; its $K$-th row contains four row vectors of three components each
\[
\left[ \begin{array}{c|c|c|c} \mathbf n_1^\top & \mathbf n_2^\top & \mathbf n_3^\top & \mathbf n_4^\top\end{array}\right]
\]
so that $\mathbf n_\ell$ is the normal vector to the face $e^K_\ell$, {\em pointing outwards} and such that $|\mathbf n_\ell|=|e^K_\ell|$.
\end{itemize}

\section{Volume integrals}

\paragraph{Pseudo-matlab notation.} In order to exploit the vectorization capabilities of Matlab, we will use the following notation to describe some particular operations. First of all, for a function $f(x,y,z)$, we will automatically assume that it is vectorized and can thus be simultaneously evaluated in many points stored in equally sizes matrix, so that $f(\mathrm X,\mathrm Y,\mathrm Z)$ is a matrix with the same size as $\mathrm X$, $\mathrm Y$, and $\mathrm Z$. Our default will be that vectors are column vectors. Whenever the sizes of a column vector $\mathbf u$ and a matrix $\mathrm A$ are compatible, we will write $\mathbf u^\top \odot \mathrm A:=\mathrm A \, \mathrm{diag}(\mathbf u), $ and $ \mathbf u\odot \mathrm A:=\mathrm{diag}(\mathbf u)\,\mathrm A$. At the entry level, these are the operations
\[
(\mathbf u^\top\odot\mathrm A)_{ij}=u_j\mathrm A_{ij}\qquad\mbox{and}\qquad (\mathbf u\odot\mathrm A)_{ij}=u_i\mathrm A_{ij},
\]
which can be easily performed using Matlab's {\tt bsxfun} utility. Also, the Kronecker product will be used in the following particular situation:
\[
\mathbf c^\top\otimes \mathrm A=\left[\begin{array}{c|c|c|c} c_1\mathrm A & c_2\mathrm A &\cdots & c_N \mathrm A\end{array}\right].
\]
Finally, given a matrix $\mathrm A$, $\mathbf a_i^\top:=\mathrm{row}(\mathrm A,i)$ will be used to denote the row vector corresponding to the $i$-th row of $\mathrm A$.

\paragraph{Three dimensional quadrature.} To compute or approximate element integrals we will consider a {\em quadrature formula} on the reference element $\widehat K$:
\begin{equation}\label{eq:3.1}
\int_{\widehat K}\widehat \phi \approx \smallfrac16\sum_{q=1}^\Nnd \widehat\omega_q\widehat\phi(\widehat{\mathbf p}_q),
\qquad\widehat{\mathbf p}_q:=(\widehat x_q,\widehat y_q,\widehat z_q), \qquad \sum_{q=1}^\Nnd \widehat \omega_q = 1.
\end{equation}
Such formulas can be easily found in the literature \cite{Felippa:2004, ZhCuLi:2009}. We will find it convenient to store the quadrature points with their barycentric coordinates $(1-\widehat x_q-\widehat y_q-\widehat z_q,\widehat x_q,\widehat y_q,\widehat z_q)$ in a $\Nnd\times 4$ matrix $\Lambda$ (rows correspond to quadrature points). To integrate on a general element we use the mapping $F_K:\widehat K\to K$ of \eqref{eq:2.0} and proceed as follows
\begin{equation}\label{eq:3.2}
\int_K\phi =\mathrm{det}\,\mathrm B_K\int_{\widehat K}\phi\circ F_K \approx |K| \sum_{q=1}^\Nnd\widehat\omega_q\phi(\mathbf p_q^K), \qquad \mathbf p_q^K:= F_K(\widehat{\mathbf p}_q).
\end{equation}

\paragraph{Piecewise polynomials.} 
The bases of the local polynomial spaces $\mathcal P_3(K)$, containing $d_3=d_3(k):={k+3\choose 3}$ elements, will be obtained by pushing forward a basis in the reference element. For this we will use the three dimensional Dubiner basis \cite{Kirby:2010}, which is given in the enlarged element $2\widehat K-(1,1,1)^\top$, where it is $L^2$-orthogonal. The Dubiner basis is evaluated using a Duffy-type transformation and Jacobi polynomials. What is needed for HDG is the evaluation of the Dubiner basis $\{\widecheck P_i\}$ and of its partial derivatives $\partial_\alpha \widecheck P_i$. {\em We will assume that the basis is ordered in hierarchical form}, that is, polynomials of degree $k$ are stored after all polynomials of degree $k-1$ for all $k$. We then consider the local bases $\{ P_i^K\}$ given by the relations
\begin{equation}\label{eq:3.3}
P_i^K\circ F_K = \widehat P_i:=\widecheck P_i(2\punto-(1,1,1)^\top).
\end{equation}
An element of the space $W_h$ will be usually stored as a $d_k\times \Nelt$ matrix, where each column contains the coefficients of the local polynomial function in the local basis.

\paragraph{Source terms.} We first extract all nodal information of the grid in three $4\times \Nelt$ matrices $\mathrm X^{\mathcal T}$, $\mathrm Y^{\mathcal T}$, $\mathrm Z^{\mathcal T}$. For instance, the element $\mathrm X^{\mathcal T}_{i,K}$ contains the $x$ coordinate of the $i$-th node of element $K$. If $\Lambda$ is the $\Nnd\times 4$ matrix with the barycentric coordinates of the quadrature points in $\widehat K$, then
\begin{equation}\label{eq:3.4}
\mathrm X:=\Lambda \mathrm X^{\mathcal T}, \qquad 
\mathrm Y:=\Lambda \mathrm Y^{\mathcal T},\qquad
\mathrm Z:=\Lambda \mathrm Z^{\mathcal T}
\end{equation}
are $\Nnd\times \Nelt$ matrices with the coordinates of all quadrature points. Further, let us consider the $\Nnd\times d_3$ matrix
\begin{equation}\label{eq:3.5}
\mathrm P_{qj}:=\widehat P_j(\widehat{\mathbf p}_q), \qquad q=1,\ldots,\Nnd, \qquad j=1,\ldots, d_3.
\end{equation}
The computational representation of the formula (see \eqref{eq:3.2} and \eqref{eq:3.3})
\[
\int_K f P_i^K = 6|K| \int_{\widehat K} (f\circ F_K)\, \widehat P_i  \approx |K| \sum_{q=1}^\Nnd f(\mathbf p_q^K)\widehat\omega_q \widehat P_i(\widehat{\mathbf p}_q) , \qquad i=1,\ldots,d_3, \quad K\in \mathcal T_h,
\]
is given by (see \eqref{eq:3.4} and \eqref{eq:3.5})
\[
\mathbf{vol}^\top \odot \big( (\widehat{\boldsymbol\omega}\odot \mathrm P)^\top \, f(\mathrm X,\mathrm Y,\mathrm Z)\big),
\]
where $\mathbf{vol}$ is the column vector containing the volumes of all elements and $\widehat{\boldsymbol\omega}$ is a column vector with the weights of the quadrature rule.

\paragraph{Mass matrices.} In order to compute mass matrices with variable density function $m$, we use \eqref{eq:3.2}-\eqref{eq:3.3} and write
\[
\int_K m\, P_i^K\,P_j^K\approx |K| \sum_{q=1}^\Nnd m(\mathbf p_q^K) \Big( \widehat \omega_q \widehat P_i(\widehat{\mathbf p}_q)
\widehat P_j(\widehat{\mathbf p}_q)\Big), \qquad i,j,=1,\ldots, d_3, \quad K\in \mathcal T_h.
\]
We then evaluate the density at all quadrature points \eqref{eq:3.4} to get an $\Nnd\times \Nelt$ matrix, already weighted by the element volumes, 
\begin{subequations}\label{eq:3.55}
\begin{equation}
\mathrm M=\mathbf{vol}^\top \odot m(\mathrm X,\mathrm Y,\mathrm Z),
\end{equation}
and finally loop over quadrature nodes using the rows of \eqref{eq:3.5}
\begin{equation}
\sum_{q=1}^\Nnd \mathbf m_q^\top \otimes (\widehat\omega_q\mathbf p_q\,\mathbf p_q^\top), \qquad  \mathbf m_q^\top=\mathrm{row}(\mathrm M,q),\quad \mathbf p_q^\top=\mathrm{row}(\mathrm P,q).
\end{equation}
\end{subequations}
 The result comes out as a $d_3\times (d_3\Nelt)$ matrix that can be easily reshaped to a $d_3\times d_3\times \Nelt$ array.

\paragraph{Convection matrices.} We start by computing three $d_3\times d_3$ matrices in the reference element
\begin{equation}\label{eq:3.6}
\widehat{\mathrm C}^\star_{ij} :=\int_{\widehat K} \widehat P_i \partial_{\widehat\star}\, \widehat P_j =\smallfrac16\sum_{q=1}^\Nnd \widehat\omega_q \widehat P_i(\widehat{\mathbf p}_q) \partial_{\widehat\star} \widehat P_j(\widehat{\mathbf p}_q), \qquad i,j=1,\ldots,d_3, \quad \star\in \{x,y,z\}.
\end{equation}
(We assume that the quadrature rule is of sufficiently high order to compute these matrices exactly.)
To do this, we require the matrix $\mathrm P$ in \eqref{eq:3.5} plus three matrices with derivatives of the basis functions in the reference element
\begin{equation}\label{eq:3.66}
\mathrm P_{qi}^\star:=(\partial_{\widehat\star} \widehat P_i)(\widehat{\mathbf p}_q) \qquad q=1.\ldots,\Nnd, \quad i=1,\ldots, d_3, \quad \star\in \{x,y,z\}. 
\end{equation}
Then, 
\begin{equation}\label{eq:3.7}
\widehat{\mathrm C}^\star=\smallfrac16 \mathrm P^\top \,\mathrm{diag}(\widehat{\boldsymbol\omega})\,\mathrm P^\star=
\smallfrac16 (\widehat{\boldsymbol\omega}\odot\mathrm P)^\top \mathrm P^\star.
\end{equation}
Next, we deal with the elements of the associated Piola transform. The elements of the $3\times 3$ matrices
\begin{equation}\label{eq:3.8}
\mathrm{det}\,\mathrm B_K\mathrm B_K^{-\top}=\left[ \begin{array}{ccc} a_{xx}^K & a_{xy}^K & a_{xz}^K \\[1.5ex] a_{yx}^K & a_{yy}^K & a_{yz}^K \\[1.5ex] a_{zx}^K & a_{zy}^K & a_{zz}^K\end{array}\right],\qquad K\in \mathcal T_h,
\end{equation}
can be computed using the coordinates of the vertices counted by elements (these are the rows of the matrices $\mathrm X^{\mathcal T}$, $\mathrm Y^{\mathcal T}$, and $\mathrm Z^{\mathcal T}$) using the formulas:
\begin{eqnarray*}
a_{xx} &=& (y_3-y_1) (z_4-z_1)-(y_4-y_1) (z_3-z_1),\\
a_{xy} &=& (y_4-y_1) (z_2-z_1)-(y_2-y_1) (z_4-z_1),\\
a_{xz} &=& (y_2-y_1) (z_3-z_1)-(y_3-y_1) (z_2-z_1),\\
a_{yx} &=& (x_4-x_1) (z_3-z_1)-(x_3-x_1) (z_4-z_1),\\
a_{yy} &=& (x_2-x_1) (z_4-z_1)-(x_4-x_1) (z_2-z_1),\\
a_{yz} &=& (x_3-x_1) (z_2-z_1)-(x_2-x_1) (z_3-z_1),\\
a_{zx} &=& (x_3-x_1) (y_4-y_1)-(x_4-x_1) (y_3-y_1),\\
a_{zy} &=& (x_4-x_1) (y_2-y_1)-(x_2-x_1) (y_4-y_1),\\
a_{zz} &=& (x_2-x_1) (y_3-y_1)-(x_3-x_1) (y_2-y_1).
\end{eqnarray*}
(Reference to $K$ has been dropped to simplify the expression.) A simple change of variables leads to
\[
\int_K P_i^K \partial_\star P_j^K = \sum_{\#\in \{x,y,z\}} a_{\star\#}^K \int_{\widehat K} \widehat P_i\partial_{\widehat\#}\widehat P_j= \sum_{\#\in \{x,y,z\}} a_{\star\#}^K \widehat{\mathrm C}_{ij}^\#,
\]
which, using the matrices \eqref{eq:3.7}, can be implemented with Kronecker products
\begin{equation}\label{eq:3.56}
\sum_{\#\in \{x,y,z\}} \mathbf a_{\star\#}^\top \otimes \widehat{\mathrm C}^\#,\qquad \star\in \{x,y,z\}.
\end{equation}
The result are three $d_3\times (d_3\Nelt) \equiv d_3\times d_3\times \Nelt$ matrices.


\section{Surface integrals}\label{sec:4}

\paragraph{Integrals on faces.}
Two dimensional quadrature rules will be given in the reference element $\widehat K_2$, using points and weights so that
\[
\int_{\widehat K_2} \widehat\phi \approx \smallfrac12 \sum_{r=1}^\Nqd \varpi_r \widehat\phi(\widehat{\mathbf q}_r),
\qquad \widehat{\mathbf q}_r=(\widehat s_r,\widehat t_r), \qquad 
\sum_{r=1}^\Nqd\varpi_r=1.
\]
To compute an integral on $e\in \mathcal E_h$, we simply parametrize from $\widehat K_2$ and proceed accordingly:
\[
\int_e \phi = 2 |e| \int_{\widehat K_2} \phi\circ\boldsymbol\phi_e\approx |e|\sum_{r=1}^\Nqd \varpi_r \phi(\mathbf q_r^e) \qquad\mbox{with}\qquad \mathbf q_r^e:=\boldsymbol\phi_e(\widehat{\mathbf q}_r).
\]
For practical purposes, we will keep the barycentric coordinates of the quadrature points $(1-s_r-t_r,s_r,t_r)$ in an $\Nqd\times 3$ matrix $\Xi$.

\paragraph{Integrals on boundaries of tetrahedra.}

In many cases we will be integrating on a face that is given with geometric information of an adjacent tetrahedron. The quadrature points $\widehat{\mathbf q}_r$ lead to four groups of quadrature points on the faces of $\widehat K$ (see \eqref{eq:2.5}):
\[
\begin{array}{l}
\widehat{\mathbf q}_r^1 :=(\widehat s_r,\widehat t_r,0),\\[1.5ex]
\widehat{\mathbf q}_r^2 :=(\widehat s_r,0,\widehat t_r),\\[1.5ex]
\widehat{\mathbf q}_r^3 :=(0,\widehat s_r,\widehat t_r),\\[1.5ex]
\widehat{\mathbf q}_r^4 :=(\widehat s_r,\widehat t_r,1-\widehat s_r-\widehat t_r),
\end{array}
\]
For a given $\psi:K\to \mathbb R$, we can approximate (see \eqref{eq:2.2} and \eqref{eq:2.3})
\begin{equation}\label{eq:4.1}
\int_{e^K_\ell}\psi \approx |e^K_\ell| \sum_{r=1}^\Nqd \varpi_r\psi (\mathbf q_{r,\ell}^K),  \mbox{ with }  \mathbf q_{r,\ell}^K := \mathrm F_K(\widehat{\mathbf q}_r^\ell)=\boldsymbol\phi_{e^K_\ell}(  F_{\mathrm{perm}(K,\ell)}(\widehat{\mathbf q}_r)),
\end{equation}
and thus
\begin{equation}\label{eq:4.2}
\int_{\partial K}\psi \approx \sum_{\ell=1}^ 4 |e^K_\ell| \sum_{r=1}^\Nqd \varpi_r\psi (\mathbf q_{r,\ell}^K) .
\end{equation}
Note that the use of the permutation index $\mathrm{perm}(K,\ell)$ in \eqref{eq:4.1}
factors out the natural parametrization of $e=e_\ell^K$ on the left, which will be necessary for functions on $e$ that are defined by pushing forward functions on $\widehat K_2$, that is, for functions $\psi$ such that we can evaluate $\psi\circ\boldsymbol\phi_e$.

\paragraph{Bases, faces, and boundaries.} Our starting point is the Dubiner basis $\{ \widecheck D_i\}$ \cite{Kirby:2010}, which is orthogonal in the enlarged element $2\widehat K_2-(1,1)^\top$. We assume it to be given in hierarchical form. The elements of $M_h$ will be described via their coefficients in the basis $D_i^e$, where
\[
D_i^e\circ\boldsymbol\phi_e=\widehat D_i=\widecheck D_i(2\punto-(1,1)^\top), \qquad i=1,\ldots,d_2, \quad e\in \mathcal E_h,\quad  d_2= d_2(k):={k+2\choose 2},
\]
so that they are stored in form of a $d_2\times \Nfc$ matrix (recall that $\Nfc=\# \mathcal E_h$).

\paragraph{Types of boundary integrals.} There will be three different kinds of integrals on $\partial K$: (a) products of traces of polynomials on $K$, (b) products of piecewise polynomials defined on $\partial K$, (c) products of traces of polynomials of $K$ by piecewise polynomials on $\partial K$. Each of these integrals will involve some kind of piecewise constant weight function. Piecewise constant functions on the boundaries of the elements (with different values on internal faces) will be described with $4\times \Nfc$ matrices. We will be using four examples of this kind of functions:
\begin{equation}\label{eq:4.22}
\begin{array}{r}
\mathrm T^K_\ell:=\tau_\ell^K |e^K_\ell|, \\ \mathrm N_{\ell,K}^\star := n_{\ell,\star}^K=\nu_{\ell,\star}^K|e^K_\ell|,
\end{array} \quad \ell=1,\ldots,4, \quad K\in \mathcal T_h, \quad \star\in \{x,y,z\}.
\end{equation}
Here $\mathbf n_\ell^K=(n_{\ell,x}^K,n_{\ell,y}^K,n_{\ell,z}^K)$ is the normal vector on the $\ell$-th face of $K$, with the normalization $|\mathbf n_\ell^K|=|e^K_\ell|$ (see Section \ref{sec:2}). The information of these four piecewise constant functions is readily available in the enhanced geometric data structure.

\paragraph{Type (a) matrices.} We can compute the integrals
\begin{equation}\label{eq:4.3}
\int_{\partial K} \tau_K P_i^K P_j^K = \sum_{\ell=1}^4 \mathrm T^K_\ell \Big( \sum_{r=1}^\Nqd \widehat P_i(\widehat{\mathbf q}_r^\ell) \varpi_r \widehat P_j(\widehat{\mathbf q}_r^\ell)\Big) \qquad i,j=1,\ldots,d_3, \quad K\in \mathcal T_h,
\end{equation}
using a sufficiently precise quadrature rule. If we consider the matrices
\begin{equation}\label{eq:4.4}
\mathrm P_{ri}^\ell:=\widehat P_i(\widehat{\mathbf q}_r^\ell), \qquad r=1,\ldots,\Nqd, \quad i=1,\ldots,d_3, \qquad \ell\in \{1,2,3,4\},
\end{equation}
then \eqref{eq:4.3} can be computed as
\begin{equation}\label{eq:4.5}
\sum_{\ell=1}^4 \mathbf t_\ell^\top \otimes \Big( (\boldsymbol\varpi\odot \mathrm P^\ell)^\top \mathrm P^\ell\Big), \qquad \mathbf t_\ell^\top=\mathrm{row}(\mathrm T,\ell)
\end{equation}
using the matrix $\mathrm T$ in \eqref{eq:4.22}.  The result is a $d_3\times (d_3\Nelt)\equiv d_3\times d_3\times \Nelt$ array.

\paragraph{Type (b) matrices.} To compute the matrices
\[
\tau_\ell^K\int_{e^K_\ell}  D_i^{e^K_\ell} D_j^{e^K_\ell}= \mathrm T_\ell^K \sum_{r=1}^\Nqd \varpi_r \widehat D_i(\widehat{\mathbf q}_r) \widehat D_j(\widehat{\mathbf q}_r) , \qquad 
\begin{array}{l} i,j=1,\ldots,d_2, \\
 K\in \mathcal T_h, \quad \ell\in \{1,2,3,4\},
\end{array}
\]
we compute the matrices
\begin{equation}\label{eq:4.7}
\mathrm D_{ri}:= \widehat D_i(\widehat{\mathbf q}_r) \qquad r=1,\ldots,\Nqd, \quad i=1,\ldots, d_2
\end{equation}
and mix them in the form
\begin{equation}\label{eq:4.8}
\mathbf t_\ell^\top \otimes \big( (\boldsymbol\varpi\odot\mathrm D)^\top\mathrm D\big), \qquad   \mathbf t_\ell^\top=\mathrm{row}(\mathrm T,\ell).
\end{equation}
A simpler option is taking advantage of the fact that the Dubiner basis is orthogonal, so these computations yield diagonal matrices. The result are four $d_2\times (d_2\Nelt)\equiv d_2\times d_2\times \Nelt$ matrices. They will be the diagonal blocks of a $(4d_2)\times (4d_2)\times \Nelt$ matrix that will be used in the local solvers.

\paragraph{Type (c) matrices.} Let $\boldsymbol\xi$ be a piecewise constant function on the set of boundaries of the elements (in practice, one of the functions described in \eqref{eq:4.22}). Let $\boldsymbol\xi_\mu$ be the piecewise constant functions given by
\begin{equation}\label{eq:4.9}
\xi_{\ell,\mu}^K:= \xi_{\ell}^K\mathbf 1_{\mathrm{perm}(K,\ell)=\mu}, \qquad \ell\in\{1,2,3,4\}, \qquad K\in \mathcal T_h, \qquad \mu\in \{1,2,3,4,5,6\}.
\end{equation}
Following \eqref{eq:4.1},
we can then compute
\begin{eqnarray*}
\frac{\xi_\ell^K}{|e^K_\ell|}\int_{e^K_\ell} D_i^{e^K_\ell} P_j^K &=&\xi_\ell^K \sum_{r=1}^\Nqd \widehat D_i(F_{\mathrm{perm}(K,\ell)}(\widehat{\mathbf q}_r))\varpi_r \widehat P_j(\widehat{\mathbf q}_r^\ell)\\
&=& \sum_{\mu=1}^6 \xi_{\ell,\mu}^K\Big(\sum_{r=1}^\Nqd \widehat D_i(F_\mu(\widehat{\mathbf q}_r)) \,\varpi_r \,\widehat P_j(\widehat{\mathbf q}_r^\ell)\Big),
\end{eqnarray*}
for $i=1,\ldots,d_2$, $j=1,\ldots,d_3$, $K\in \mathcal T_h$ and $\ell\in \{1,2,3,4\}$. Using \eqref{eq:4.4}, the matrices
\[
\mathrm D_{ri}^\mu := \widehat D_i(F_\mu(\widehat{\mathbf q}_r)) \qquad r=1,\ldots,\Nqd, \quad i=1,
\ldots,d_2, \quad \mu\in \{1,\ldots,6\},
\]
and \eqref{eq:4.9}, the previous computation reduces to
\begin{equation}\label{eq:4.10}
\sum_{\mu=1}^6 \boldsymbol\xi_{\ell,\mu}^\top \otimes ((\boldsymbol\varpi\odot \mathrm D^\mu)^\top \mathrm P^\ell), \qquad \boldsymbol\xi_{\ell,\mu}^\top=\mathrm{row}(\boldsymbol\xi_\mu,\ell) \qquad \ell=1,2,3,4.
\end{equation}
The result is four $d_2\times (d_3\Nelt)\equiv d_2\times d_3\times \Nelt$ matrices that are stored as a single $(4d_2)\times d_3\times \Nelt$ array, by stacking the blocks for $\ell=1,2,3,4$ on top of each other ($\ell=1$ on top).

\section{Local solvers}\label{sec:5}

The local solvers that we next define are related to the pair of discrete equations \eqref{eq:1.2a}-\eqref{eq:1.2b}.

\paragraph{Matrices and bilinear forms.}
In order to recognize the matrices that we have computed with terms in the bilinear forms of the HDG method, we need some notation. We consider the space
\[
\mathcal R_k(\partial K):=\prod_{e\in \mathcal E(K)} \mathcal P_k(e), \qquad \mathrm{dim}\,\mathcal R_k(\partial K)=4 d_2.
\]
The degrees of freedom for this last space are organized by taking one face at a time in the order they are given by {\tt T.facebyele}. For (non-symmetric) bilinear forms we will use the convention that the bilinear form $b(u,v)$ is related to the matrix $b(U_j,V_i)$, where $\{U_j\}$ is a basis of the space of $u$ and $\{V_i\}$ is a basis of the space for $v$. This is equivalent to saying that the unknown will always be placed as the left-most argument in the bilinear form and the test function will occupy the right-most location.

\paragraph{Volume terms.} We start by computing mass matrices associated to two functions ($\kappa^{-1}$ and $c$), and the three convection matrices:
\[
\mathrm M_{\kappa^{-1}}^K, \qquad \mathrm M_c^K, \qquad \mathrm C_x^K, \qquad \mathrm C_y^K, \qquad \mathrm C_z^K,
\]
where (see \eqref{eq:3.55} and \eqref{eq:3.56})
\[
(\mathrm M_m^K)_{ij}=\int_K m P_i^K P_j^K, \qquad (\mathrm C_\star^K)_{ij}=\int_K P_i^K \partial_{\star} P_j^K.
\]
Each of these matrices is $d_3\times d_3\times \Nelt$. They correspond to the bilinear forms
\[
(m\,u_h,v_h)_K, \qquad (\partial_{\star} u_h,v_h)_K, \qquad u_h, v_h \in \mathcal P_k(K).
\]

\paragraph{Surface terms.} We next compute all matrices related to integrals on interfaces:
\[
\tau\mathrm{PP}^K, \qquad \tau\mathrm{DP}^K,  \qquad n_x\mathrm{DP}^K, \qquad n_y\mathrm{DP}^K, \qquad n_z\mathrm{DP}^K,\qquad \tau\mathrm{DD}^K.
\]
The first of these arrays is $d_3\times d_3\times \Nelt$, the next four are $4d_2\times d_3\times \Nelt$ and the last one is $4d_2\times 4d_2\times \Nelt$. The first matrix and associated bilinear form (see \eqref{eq:4.5}) are
\[
\tau\mathrm{PP}^K_{ij} =\int_{\partial K} \tau_K\, P_i^K P_j^K, \qquad \langle \tau u_h,v_h\rangle_{\partial K}, \qquad u_h,v_h \in \mathcal P_k(K).
\]
The second one (see \eqref{eq:4.10}) corresponds to the bilinear form
\[
\langle \tau_K u_h,\widehat v_h\rangle_{\partial K}, \qquad u_h\in \mathcal P_k(K), \quad \widehat v_h \in \mathcal R_k(\partial K),
\]
or equivalently to $\langle \tau_K u_h,\widehat v_h\rangle_e$, for $u_h\in \mathcal P_k(K)$, $\widehat v_h \in \mathcal P_k(e)$, and $e \in \mathcal E(K).$
The matrices associated to the components of the normal vector $\boldsymbol\nu=(\nu_x,\nu_y,\nu_z)$ (see \eqref{eq:4.10} again) are related to the bilinear forms
\[
\langle \nu_\star u_h,\widehat v_h\rangle_{\partial K}, \qquad u_h\in \mathcal P_k(K), \quad \widehat v_h \in \mathcal R_k(\partial K), \qquad \star\in \{x,y,z\}.
\]
The last matrix (see \eqref{eq:4.8}) corresponds to
\[
\langle \tau \widehat u_h,\widehat v_h\rangle_{\partial K}, \qquad \widehat u_h,\widehat v_h \in \mathcal R_k(\partial K),
\]
and is therefore block diagonal. Finally we compute the vectors of tests of $f$ with the basis elements of $\mathcal P_k(K)$: $\mathbf f^K \in \mathbb R^{d_3}$.

\paragraph{Matrices related to local solvers.}
The $4d_3\times 4d_3\times \Nelt$ array and the $4d_3\times 4d_2\times \Nelt$ array with respective slices
\begin{equation}\label{eq:5.1}
\mathbb A_1^K:=\left[\begin{array}{cccc} \mathrm M_{\kappa^{-1}}^K & \mathrm O &\mathrm O & -(\mathrm C_x^K)^\top \\ 
 \mathrm O & \mathrm M_{\kappa^{-1}}^K & \mathrm O & -(\mathrm C_y^K)^\top\\
\mathrm O & \mathrm O & \mathrm M_{\kappa^{-1}}^K & -(\mathrm C_z^K)^\top\\
\mathrm C_x^K & \mathrm C_y^K & \mathrm C_y^K & \mathrm M_c^K+\tau\mathrm{PP}^K
\end{array}\right],\qquad  \mathbb A_2^K:=\left[\begin{array}{c} (n_x\mathrm{DP}^K)^\top\\ (n_y\mathrm{DP}^K)^\top \\
(n_z\mathrm{DP}^K)^\top \\
-(\tau\mathrm{DP}^K)^\top\end{array}\right],
\end{equation}
are the matrix representations of the bilinear forms
\begin{eqnarray*}
a_1^K : \big(\mathcal P_k(K)^3\times \mathcal P_k(K)\big)\times\big(\mathcal P_k(K)^3\times \mathcal P_k(K)\big)&\longrightarrow & \mathbb R,\\
a_2^K :  \mathcal R_k(\partial K)\times \big(\mathcal P_k(K)^3\times \mathcal P_k(K)\big)&\longrightarrow & \mathbb R,
\end{eqnarray*}
given by
\begin{eqnarray*}
a_1^K((\boldsymbol q_h,u_h),(\boldsymbol r_h,w_h)) &:=& (\kappa^{-1}\boldsymbol _h,\boldsymbol r_h)_K -(u_h,\nabla\cdot\boldsymbol r_h)_K \\
& & + (\nabla\cdot\boldsymbol q_h,w_h)_K+(c\,u_h,w_h)_K+\langle \tau u_h,v_h\rangle_{\partial K},\\
a_2^K(\widehat u_h, (\boldsymbol r_h,w_h))&:=& \langle \widehat u_h,\boldsymbol r_h \cdot\boldsymbol\nu\rangle_{\partial K}-\langle \tau \widehat u_h,w_h\rangle_{\partial K}.
\end{eqnarray*}
We also consider the $4d_3\times \Nelt$ matrix with columns
\begin{equation}\label{eq:5.2}
\mathbb A_f^K:=\left[\begin{array}{ccc} \mathbf 0 \\ \mathbf 0 \\ \mathbf 0\\ \mathbf f^K\end{array}\right],
\end{equation}
If $\widehat u_h \in M_h$ is known, we can solve the local problems looking for $\boldsymbol q_h \in \boldsymbol V_h= W_h^3$ and $u_h \in W_h$, satisfying \eqref{eq:1.2a}-\eqref{eq:1.2b}. Representing
 $\widehat u_h|_{\partial K}\in \mathcal R_k(\partial K)$  with a vector $\mathbf u_{\partial K} \in \mathbb R^{4d_2}$,  the matrix representation of this local solution is
\begin{equation}\label{eq:5.3}
\left[\begin{array}{c} \mathbf q_K \\ \mathbf u_K\end{array}\right]=-(\mathbb A_1^K)^{-1} \mathbb A_2^K \mathbf u_{\partial K} + (\mathbb A_1^K)^{-1}\mathbb A_f^K \in \mathbb R^{4d_3}.
\end{equation}
Note that once the local matrices have been computed, the construction of the three dimentional arrays \eqref{eq:5.1} can be easily carried out by stacking the already created three dimensional arrays.

\paragraph{Flux operators.}
Consider now the $4d_2\times 4d_3\times \Nelt$ array with slices
\begin{equation}\label{eq:5.4}
\mathbb A_3^K:=\left[\begin{array}{cccc} n_x\mathrm{DP}^K & n_y\mathrm{DP}^K & n_z\mathrm{DP}^K &
\tau\mathrm{DP}^K\end{array}\right],
\end{equation}
the $4d_2\times 4d_2\times \Nelt$ array with slices
\begin{equation}\label{eq:5.5}
\mathbb C^K:= \mathbb A_3^K ( \mathbb A_1^K)^{-1} \mathbb A_2^K+\tau\mathrm{DD}^K,
\end{equation}
and the $4d_2\times \Nelt$ matrix with columns
\begin{equation}\label{eq:5.6}
 \mathbb C_f^K:=\mathbb A_3^K (\mathbb A_1^K)^{-1}\mathbb A_f^K.
\end{equation}
The meaning of these matrices  can be made clear by looking at boundary fluxes. 
Given $(\boldsymbol q_h,u_h,\widehat u_h) \in \boldsymbol V_h\times W_h \times M_h$ --satisfying  equations \eqref{eq:1.2a} and \eqref{eq:1.2b}--, the HDG method is based on the construction of the flux function
\[
\Phi_K:=-\boldsymbol q_h\cdot\boldsymbol\nu-\tau (u_h-\widehat u_h): \partial K \to \mathbb R.
\]
Instead of this quantity, we pay attention to how it creates a linear form
\[
\mathcal R_k(\partial K) \ni \widehat v_h \longmapsto -\langle \boldsymbol q_h\cdot\boldsymbol\nu+\tau (u_h-\widehat u_h),\widehat v_h\rangle_{\partial K}=-\langle \boldsymbol q_h\cdot\boldsymbol\nu+\tau u_h,\widehat v_h\rangle_{\partial K}+\langle \tau \widehat u_h,\widehat v_h\rangle_{\partial K},
\]
whose matrix representation is
\begin{eqnarray}\nonumber
-\mathbb A_3^K \left[\begin{array}{c} \mathbf q_K \\ \mathbf u_K\end{array}\right]+\tau\mathrm{DD}^K\mathbf u_{\partial K} &=&
\mathbb A_3^K(\mathbb A_1^K)^{-1} \mathbb A_2^K \mathbf u_{\partial K}-\mathbb A_3^K(\mathbb A_1^K)^{-1}\mathbb A_f^K+\tau\mathrm{DD}^K\mathbf u_{\partial K}\\
\label{eq:5.7}
&=&
\mathbb C^K \mathbf u_{\partial K}-\mathbb C_f^K,
\end{eqnarray}
where $\mathbf u_{\partial K}$ is the vector of degrees of freedom of $\widehat u_h|_{\partial K}$.

\paragraph{Note on implementation.} Construction of the local solvers \eqref{eq:5.5} and \eqref{eq:5.6}, as well as recovery of internal values using \eqref{eq:5.3}, requires looping over elements. However, this can be easily done in parallel, since at this stage there is no interconnection between elements. Note that we have avoided looping over elements in all previous computations, requiring frequent access to coefficients and geometric features.

\section{Boundary conditions and global solver}

\paragraph{Dirichlet boundary conditions.} The discrete Dirichlet boundary conditions require finding
the decompositions
\[ 
\widehat u_h|_e=\sum_{j=1}^{d_2} u_j^e D_j^e, \qquad e\in \mathcal E_h^D
\]
by solving the system
\[
\sum_j \left(\int_e D_i^e D_j^e \right)u_{j}^e=\int_e D_i^e u_D, \qquad  i=1,\ldots,d_2, \qquad e\in \mathcal E_h^{D}.
\]
Using a quadrature rule on the reference element, and parametrizing from it, we have to solve the approximate system
\begin{equation}\label{eq:6.1}
|e| \sum_{j=1}^{d_2} \Big(\sum_{r=1}^{\Nqd} \varpi_r \widehat D_i(\widehat{\mathbf q}_r)\widehat D_j(\widehat{\mathbf q}_j)\Big) u_j^e=|e| \sum_{r=1}^{\Nqd} \varpi_r \widehat D_i(\widehat{\mathbf q}_r) u_D(\mathbf q_r^e).
\end{equation}
Evaluation of the data function $u_D$ at all the quadrature points is done with a similar strategy to the one used for source terms \eqref{eq:3.4}. We start by organizing nodal information in three $3\times \Ndir$ (recall that $\Ndir=\# \mathcal E_h^D$) matrices $\mathrm X^D$, $\mathrm Y^D$, $\mathrm Z^D$, each of the containing the corresponding coordinates of the nodes of each of the Dirichlet faces. If $\Xi$ is the $\Nqd\times 3$ matrix with the barycentric coordinates of the quadrature points in $\widehat K_2$ (Section \ref{sec:4}), then the $\Nqd\times \Ndir$ matrices
\begin{equation}\label{eq:6.2}
\mathrm X^{\mathrm{dir}}:=\Xi\,\mathrm X^D, \qquad\mathrm Y^{\mathrm{dir}}:=\Xi\,\mathrm Y^D,\qquad\mathrm Z^{\mathrm{dir}}:=\Xi\,\mathrm Z^D,
\end{equation}
contain the coordinates of the quadrature points on the Dirichlet faces. Using the matrix in \eqref{eq:4.7} (see also \eqref{eq:4.8}),  it is clear that \eqref{eq:6.1} can be implemented by solving a system with multiple right-hand sides
\begin{equation}\label{eq:6.3}
\widehat{\mathbf u}_D:= \big( (\boldsymbol\varpi\odot\mathrm D)^\top\mathrm D\big)^{-1} (\boldsymbol\varpi\odot\mathrm D)^\top u_D(\mathrm X^{\mathrm{dir}},\mathrm Y^{\mathrm{dir}},\mathrm Z^{\mathrm{dir}}).
\end{equation}
The result is a $d_2\times \Ndir$ matrix.

\paragraph{Neumann boundary conditions.} As opposed to Dirichlet conditions (that are essential in this formulation), Neumann boundary conditions will appear in the right-hand side of the global system. Our goal is to compute the integrals (recall that $|\mathbf n_e|=|e|$)
\[
\int_e (\boldsymbol g_N\cdot\boldsymbol\nu_e) D_i^e \approx \sum_{r=1}^\Nqd \varpi_r\widehat D_i(\widehat{\mathbf q}_r)\, \boldsymbol g_N(\mathbf q_r^e) \cdot\mathbf n_e \qquad i=1,\ldots, d_2, \qquad e\in \mathcal E_h^N.
\]
If we consider matrices $\mathrm X^{\mathrm{neu}}$, $\mathrm Y^{\mathrm{neu}}$, $\mathrm Z^{\mathrm{neu}}$, defined as in \eqref{eq:6.2} (but using nodal information for Neumann faces), and if $\mathbf n_x$, $\mathbf n_y$, $\mathbf n_z$ are $\Nneu\times 1$ (recall that $\Nneu=\#\mathcal E_h^N$) column vectors with the components of the vectors $\mathbf n_e$ for $e\in \mathcal E_h^N$, then everything is done with the simple computation
\begin{equation}\label{eq:6.4}
\boldsymbol \Phi_N:=\sum_{\star\in \{x,y,z\}} \mathbf n_\star^\top \odot \big( (\boldsymbol\varpi\odot\mathrm D)^\top g_\star(\mathrm X^{\mathrm{neu}},\mathrm Y^{\mathrm{neu}},\mathrm Z^{\mathrm{neu}})\big),
\end{equation}
where $\boldsymbol g_N=(g_x,g_y,g_z)$. Note that if the Neumann boundary condition is given in a more standard way $-\mathbf q\cdot\boldsymbol\nu=g_N$, then the computation is slightly simpler
\[
\boldsymbol\Phi_N:=\mathbf{area}_N^\top \odot\big( (\boldsymbol\varpi\odot\mathrm D)^\top g_N(\mathrm X^{\mathrm{neu}},\mathrm Y^{\mathrm{neu}},\mathrm Z^{\mathrm{neu}})\big)
\]
where $\mathbf{area}_N$ contains the areas of all Neumann faces.

\paragraph{Assembly process.}
The local solvers produce a $4d_2\times 4d_2\times \Nelt$ array $\mathbb C$. We now use the {\tt sparse} MATLAB builder to assembly the global matrix. The degrees of freedom associated to face $e\in \{1,\ldots,\Nfc\}$ are
\[
\mathrm{list}(e):=(e-1)d_2+\{1,\ldots,d_2\}.
\]
The degrees of freedom associated to the faces of $K$ are thus
\[
\mathrm{dof}(K):=\{ \mathrm{list}(e_1^K),\mathrm{list}(e_2^K),\mathrm{list}(e_3^K),\mathrm{list}(e_4^K)\}.
\]
We then create two new $4d_2\times 4d_2\times \Nelt$ arrays
\[
\mathrm{Row}_{ij}^K = \mathrm{dof}(K)_i \qquad \mathrm{Col}_{ij}^K = \mathrm{dof}(K)_j, \quad\mbox{so that}\quad (\mathrm{Col}^K)^\top=\mathrm{Row}^K.
\]
The$(i,j)$  element of $\mathbb C^K$ has to be assembled at the location $(\mathrm{Row}^K_{ij},\mathrm{Col}^K_{ij})=
( \mathrm{dof}(K)_i ,  \mathrm{dof}(K)_j)$.  The result is a sparse $d_2\Nfc\times d_2\Nfc$ matrix $\mathbb H$. This matrix collects the fluxes \eqref{eq:5.7} for all the elements, with the result that opposing sign fluxes in internal faces (the normal vector points in different directions) are added.
The assembly of the {source term}, given in the matrix $\mathbb C_f$, can be carried out  using the {\tt accumarray} command. The element $(\mathbb C_f^K)_i$ has to be added to the location $\mathrm{dof}(K)_i$. The result is a vector $\mathbf F$ with $d_2\Nfc$ components. Let us consider the system at its current stage
\begin{equation}\label{eq:6.5}
\mathbb H \, \widehat{\mathbf u}=\mathbf F+ \mathbf G_N,
\end{equation}
where $\mathbf G_N$ is the $d_2\Nfc$ vector containing the elements of $\boldsymbol\Phi_N$ in the degrees of freedom corresponding to Neumann faces and zeros everywhere else. This is the matrix representation of the system \eqref{eq:1.2} with no Dirichlet boundary conditions, i.e., assuming homogeneous Neumann boundary conditions on $\Gamma_N$, the system having been written in the $\widehat u_h$ variable after local inversion of \eqref{eq:1.2a}-\eqref{eq:1.2b}.

What is left is the standard elimination of Dirichlet degrees of freedom from \eqref{eq:6.5}, namely values of Dirichlet faces are taken from \eqref{eq:6.3} and sent to the right-hand side of the system, and rows corresponding to Dirichlet degrees of freedom are ignored.

\paragraph{Reconstruction.} The solution of the resulting system is $\widehat u_h\in M_h$. Reconstruction of the other variables $(\boldsymbol q_h,u_h)$ is done by solving local problems. In matrix form, we have to solve on each $K\in \mathcal T_h$ the system
\[
\mathbb A_1^K \left[\begin{array}{c}\mathbf q_K\\ \mathbf u_K\end{array}\right]=\mathbb A^K_f-\mathbb A^K_2 \mathbf u_{\partial K}.
\]
This can be done in parallel.

\section{Add-ons}

\paragraph{Matrices for convection-diffusion problems.}
With very similar techniques, it is easy to compute convection matrices with variable coefficients
\begin{equation}\label{eq:9.1}
\int_K m\, P_i^K \partial_\star P_j^K, \qquad i,j=1,\ldots,d_3, \qquad K\in \mathcal T_h,
\end{equation}
as well as surface matrices with variable coefficients
\begin{eqnarray*}
\frac{\xi_\ell^K}{|e^K_\ell|}\int_{e^K_\ell} \alpha \, D_i^{e^K_\ell}P_j^K,& & \begin{array}{c} i=1,\ldots,d_2, \\ 
 j=1,\ldots,d_3,\end{array} \quad \ell\in \{1,2,3,4\}, \quad K\in \mathcal T_h,\\
\frac{\xi_\ell^K}{|e^K_\ell|}\int_{e^K_\ell} \alpha \, D_i^{e^K_\ell}D_j^{e^K_\ell},& & i,j,=1,\ldots,d_2, \quad \ell\in \{1,2,3,4\}, \quad K\in \mathcal T_h.
\end{eqnarray*}
These matrices are needed for  coding HDG applied to convection-diffusion problems \cite{ChCo:2012}, which needs the bilinear forms
\[
(\boldsymbol\beta\cdot\nabla u_h,w_h)_K, \quad \langle(\boldsymbol\beta\cdot\boldsymbol\nu) u_h,\widehat v_h\rangle_{\partial K},\quad \mbox{and}\quad \langle (\boldsymbol\beta\cdot\boldsymbol\nu)\widehat u_h,\widehat v_h\rangle_{\partial K}.
\]

\paragraph{Postprocessing.} If we look for $u_h^*:\Omega\to \mathbb R$ such that $u_h^*|_K\in \mathcal P_{k+1}(K)$ and for all $K\in \mathcal T_h$,
\begin{subequations}\label{eq:9.11}
\begin{alignat}{4}
(\nabla u_h^*,\nabla w_h)_K&=-( \kappa^{-1}\boldsymbol q_h,\nabla w_h)_K \qquad \forall w_h\in \mathcal P_{k+1}(K),\\
(u_h^*,1)_K &=(u_h,1)_K,
\end{alignat}
\end{subequations}
then it can be shown that this local postprocessed approximation has one additional order of convergence \cite{CoGoSa:2010}. In order to compute this postprocessing we have to use matrices of the form \eqref{eq:9.1} in the right-hand side (using an additional polynomial degree) and we need to compute local stiffness matrices
\[
\int_K \nabla P_i^K \cdot\nabla P_j^K \qquad i,j,=1,\ldots,d_3(k+1), \qquad K\in \mathcal T_h.
\]
As in the computation of the convection matrices \eqref{eq:3.56}, this can be done using geometric vectors and Kronecker products.

\paragraph{Local $L^2$ projections.}
For several different purposes, it is also convenient to have some local projections at hand. The first one is the $L^2(\Omega)$ projection on $W_h$: given $f$ we compute $f_h \in W_h$ such that
\[
(f_h,w_h)_K=(f,w_h)_K \qquad \forall w_h \in \mathcal P_k(K)\qquad \forall K\in \mathcal T_h.
\]
Using \eqref{eq:3.4} and \eqref{eq:3.5}, and up to quadrature errors, this projection is easily computed with a single instruction
\[
\big((\widehat{\boldsymbol\omega}\odot\mathrm P)^T \mathrm P\big)^{-1} (\widehat{\boldsymbol\omega}\odot\mathrm P)^\top f(\mathrm X,\mathrm Y,\mathrm Z).
\]
The $L^2$ projection on $M_h$
\begin{equation}\label{eq:9.12}
\langle f_h,\widehat v_h\rangle_e=\langle f,\widehat v_h\rangle_e \qquad \forall \widehat v_h \in \mathcal P_k(e) \qquad \forall e\in \mathcal E_h,
\end{equation}
is computed using a formula like \eqref{eq:6.3}
\[
 \big( (\boldsymbol\varpi\odot\mathrm D)^\top\mathrm D\big)^{-1} (\boldsymbol\varpi\odot\mathrm D)^\top u_D(\mathrm X^{\mathrm{all}},\mathrm Y^{\mathrm{all}},\mathrm Z^{\mathrm{all}}),
\]
where we use quadrature points on all faces of the triangulation (see \eqref{eq:6.2} for the Dirichlet case). 

\paragraph{Error functions.} Once again with very similar ideas it is easy to code the computation of errors
\[
\int_\Omega | u-u_h|^2 \qquad \sum_{e\in \mathcal E_h} |e| \int_e |\widehat u_h-u|^2,
\]
for a given function $u$ and approximations $u_h \in W_h$ and $\widehat u_h\in M_h$.

\paragraph{HDG projection.} A final projection is directly tied to the HDG method. The input is the collection $(\boldsymbol q,u)$ of a vector field an a scalar function. The output are functions $(\boldsymbol q_h,u_h)\in \boldsymbol V_h \times W_h$ satisfying
\begin{subequations}\label{eq:9.13}
\begin{alignat}{6}
(\boldsymbol q_h,\boldsymbol r)_K &=(\boldsymbol q,\boldsymbol r)_K & \qquad &\forall \boldsymbol r \in \mathcal P_{k-1}(K)^3 &\qquad & \forall K\in \mathcal T_h,\\
(u_h,w_h)_K &=(u,w_h)_K & & \forall w_h \in \mathcal P_{k-1}(K) & & \forall K\in \mathcal T_h,\\
\langle\boldsymbol q_h\cdot\boldsymbol\nu+\tau_K u_h,\widehat v_h\rangle_{\partial K} &=\langle \boldsymbol q\cdot\boldsymbol\nu+\tau\,u,\widehat v_h\rangle_{\partial K} & & \forall \widehat v_h \in \mathcal R_k(\partial K) & &\forall K\in \mathcal T_h.
\end{alignat}
\end{subequations}
It has to be understood that the first two groups of equations are void when $k=0$. If we construct a mass matrix (with constant unit mass) $\mathrm M^K$ and drop the last $d_2=\mathrm{dim}\mathcal P_k(K)-\mathrm{dim}\mathcal P_{k-1}(K)$ rows (recall that local bases are hierarchical), we obtain a $(d_3-d_2)\times d_3\times \Nelt$ matrix with slices $\widetilde{\mathrm M}^K$ 
\[
\int_K P_i^K\, P_j^K, \qquad i=1,\ldots, d_3-d_2, \qquad j=1,\ldots, d_3.
\]
Using the surface matrices of Section \ref{sec:5}, we are led to solve local linear systems with matrices:
\[
\left[\begin{array}{cccc} \widetilde{\mathrm M}^K & \mathrm O & \mathrm O & \mathrm O\\
\mathrm O & \widetilde{\mathrm M}^K & \mathrm O & \mathrm O\\
\mathrm O & \mathrm O & \widetilde{\mathrm M}^K & \mathrm O\\
\mathrm O & \mathrm O & \mathrm O & \widetilde{\mathrm M}^K\\
n_x\mathrm{DP}^K &  n_y\mathrm{DP}^K & n_z\mathrm{DP}^K & \tau\mathrm{DP}^K
\end{array}\right].
\]
The corresponding right-hand sides can be easily constructed using the techniques of previous sections.

\paragraph{BDM.} A hybridized coding of the three dimensional Brezzi-Douglas-Marini element (more properly speaking, this is an element by Brezzi-Douglas-Dur\'an-Fortin, discovered simultaneously by N\'ed\'elec) is also easily attainable. For this case, we take $k\ge 1$, define
\[
W_h:=\prod_{K\in \mathcal T_h}\mathcal P_{k-1}(K), \qquad \boldsymbol V_h:=\prod_{K\in \mathcal T_h}\mathcal P_k(K)^3
\]
and keep $M_h$ as before. The mixed BDM approximation to \eqref{eq:1.1} uses equations \eqref{eq:1.2} with two simple modifications: $\tau\equiv 0$, and equation \eqref{eq:1.2a} is only tested in $\mathcal P_{k-1}(K)^3$. At the implementation level, this means that we only need to redefine the local solvers. Since $\mathrm{dim}\mathcal P_{k-1}(K)=d_3-d_2$ and the only unknown that is in a smaller space is $u_h$, we only need to eliminate: the last $d_2$ rows and columns of $\mathbb A_1^K$, the last $d_2$ rows of $\mathbb A_2^K$ and $\mathbb A_f^K$, and the last $d_2$ columns of $\mathbb A_3^K$. (This can be done by erasing the corresponding parts of the three dimensional arrays where we have stored the HDG matrices.) All other parts of the HDG code remain untouched.

\section{Experiments}

\begin{figure}
\begin{center}
\includegraphics[width=6cm]{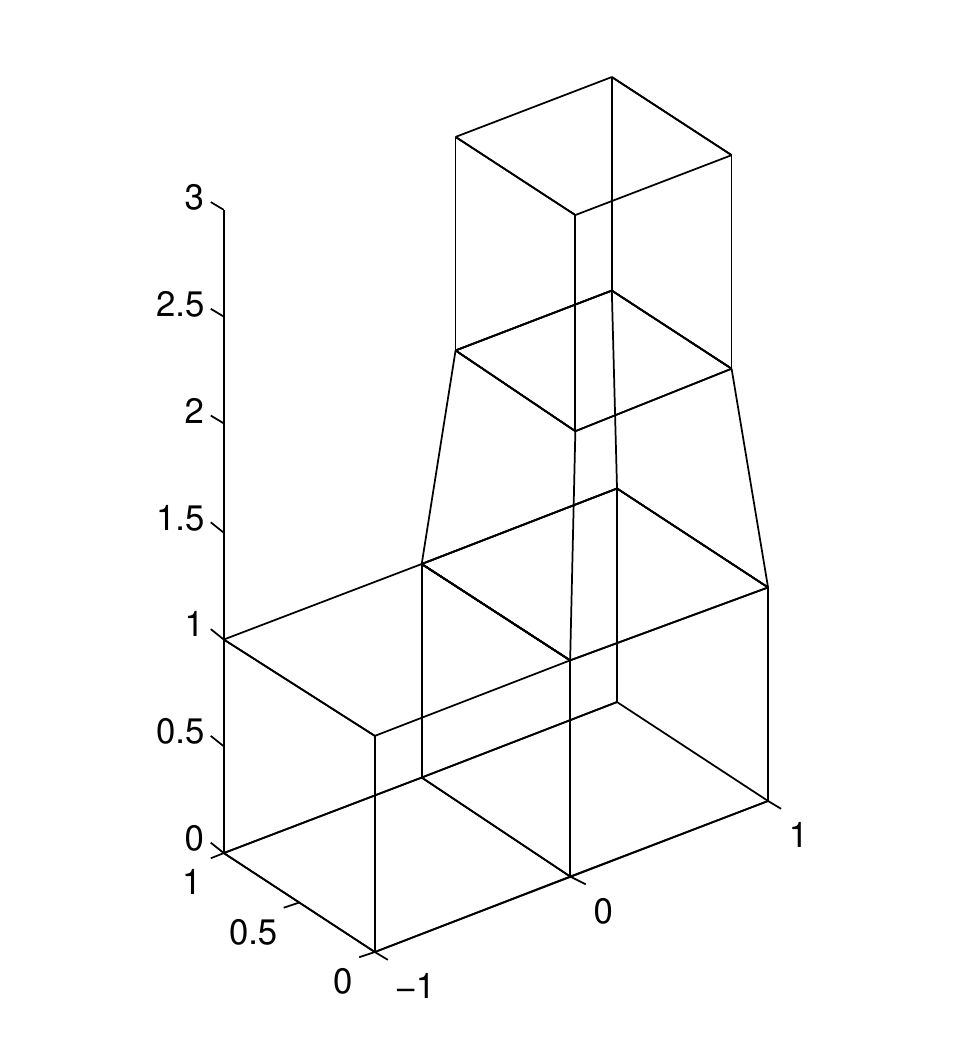}\qquad

\begin{tabular}{| c| c c c c |} 
  \hline
$\Nelt$ 	&	24 & 192 & 1536 & 12288 \\
$\Nfc$         & 66 & 456 & 3360 & 25728 \\
 \hline \end{tabular}
\end{center}
\caption{Domain for the experiments and data for the four tetrahedrizations used in the experiments.}\label{fig:1}
\end{figure}

We next give some convergence tests for the method. We take $\Omega$ to be the polyhedron sketched in Figure \ref{fig:1}. The faces of the polyhedron corresponding to $z=0$, $z=1$ and $z=3$ conform the Dirichlet boundary $\Gamma_D$. The coarsest triangulation --obtained by a tetrahedral partition of each the four hexahedra shown in Figure \ref{fig:1}--, contains $24$ elements. Three nested refinements of this partition are used. The main triangulation data are then given in Figure \ref{fig:1}. We use variable coefficients:
\[
\kappa=2+\sin x\,\sin y\,\sin z,  \qquad c=1+\smallfrac12(x^2+y^2+z^2),
\]
and take data so that $u=\sin(x\,y\,z)$ is the exact solution.

We test for several values of $k$ on the four triangulations. Tables \ref{table:1}--\ref{table:3} show relative errors
\begin{equation}\label{eq:10.1}
e^q_h:=\frac{\|\boldsymbol q-\boldsymbol q_h\|_\Omega}{\|\boldsymbol q\|_\Omega} \qquad e_h^u:=\frac{\|u-u_h\|_\Omega}{\|u\|_\Omega}, \qquad e^{\hat u}_h:=\frac{\| u-\widehat u_h\|_h}{\| u\|_h}
\end{equation}
and relative errors for superconvergent quantities
\begin{equation}\label{eq:10.2}
\varepsilon_h^u:=\frac{\| \Pi u-u_h\|_\Omega}{\| u\|_\Omega}, \qquad \varepsilon^{\hat u}_h:=\frac{\| P u-\widehat u_h\|_h}{\| u\|_h}, \qquad e^\star_h:=\frac{\| u-u_h^\star\|_\Omega}{\|u\|_\Omega}.
\end{equation}
Here: $\|u\|_h^2:=\sum_{e\in \mathcal E_h}|e| \,\| u\|_e^2$, $u_h^\star$ is the postprocessed solution defined by \eqref{eq:9.11}, $P u \in M_h$ is the $L^2(\partial\mathcal T_h)$ projection defined in \eqref{eq:9.12} and $\Pi u\in W_h$ is the scalar component of the projection $(\Pi \boldsymbol q,\Pi u)$ defined in \eqref{eq:9.13}. Theory \cite{CoGoSa:2010} shows that for smooth solutions, the errors \eqref{eq:10.1} behave like $\mathcal O(h^{k+1})$ while errors \eqref{eq:10.2} behave like $\mathcal O(h^{k+2})$ except when $k=0$, where they behave like $\mathcal O(h)$. Estimates of order of convergence for a general quantity $e_h$ are computed using the formula $\log_2 e_{h/2}/\log_2 e_h$.

\begin{table}[ht] 
\centering
   \begin{tabular}{ c c c c c c } 
  \hline\hline 
   $e_h^q$ & e.c.r. & $ e_h^u$ & e.c.r. &$e_h^{\hat u} $&e.c.r.  \\
    \hline
  6.9694e-001  &  -- & 1.2404e+000  & -- & 5.7265e-001 & -- \\ 
  4.2490e-001  & 0.71 & 8.5011e-001  & 0.55 & 3.2710e-001 & 0.81  \\ 
  2.2749e-001  &0.90 & 5.0782e-001  & 0.74 & 1.7587e-001 & 0.90 \\
  1.1739e-001  & 0.95 & 2.8336e-001  & 0.85 & 9.1779e-002 & 0.94  
   \\[1ex] \hline \end{tabular} 
\begin{tabular}{ c c c c c c } 
  \hline 
   $\varepsilon_h^u$ & e.c.r. & $ \varepsilon_h^{\hat u}$ & e.c.r. &$\varepsilon_h^\star $&e.c.r.  \\
    \hline
  1.3127e+000  &  -- & 5.4851e-001  & -- & 1.2515e+000 & -- \\ 
  6.4689e-001  & 1.02 & 2.9624e-001  & 0.89 & 8.4019e-001& 0.57  \\ 
  2.6697e-001   &1.28 & 1.5688e-001  & 0.92 & 5.0151e-001 & 0.74 \\
  1.1100e-001  & 1.27 & 8.1482e-002  & 0.95 & 2.8027e-001 & 0.84  
   \\[1ex] \hline \end{tabular} 
\caption{Errors for different triangulations (see Figure \ref{fig:1}) with the lowest order method $k=0$.}\label{table:1}
\end{table}

      \begin{table}[ht] 
\centering
   \begin{tabular}{ c c c c c c } 
  \hline\hline 
   $e_h^q$ & e.c.r. & $ e_h^u$ & e.c.r. &$e_h^{\hat u} $&e.c.r.  \\
    \hline
  1.3607e-001  &  -- & 4.2677e-001  & -- & 1.3580e-001 & -- \\ 
  3.6794e-002  & 1.89 & 1.2953e-001  & 1.72 & 3.1932e-002 & 2.09  \\ 
  9.6645e-003  &1.93 & 3.6956e-002  & 1.81 & 7.9337e-003 & 2.01 \\
  2.4878e-003  & 1.96 & 9.9288e-003  & 1.90 & 1.9876e-003 & 2.00  
   \\[1ex] \hline \end{tabular} 
\begin{tabular}{ c c c c c c } 
  \hline 
   $\varepsilon_h^u$ & e.c.r. & $ \varepsilon_h^{\hat u}$ & e.c.r. &$\varepsilon_h^\star $&e.c.r.  \\
    \hline
  2.9210e-001  &  -- & 4.6150e-002 & -- & 3.3850e-002 & -- \\ 
  4.5051e-002  & 2.70 & 5.8813e-003  & 2.97 & 4.4054e-003& 2.94  \\ 
  6.5777e-003  &2.78 & 7.6131e-004   & 2.95 & 5.5048e-004 & 3.00 \\
 8.9450e-004 & 2.88 & 9.6512e-005  & 2.98 & 6.8066e-005 & 3.02 
   \\[1ex] \hline \end{tabular} 
\caption{Errors for different triangulations (see Figure \ref{fig:1}) with the lowest order method $k=1$. All quantities in \eqref{eq:10.2} are shown to be superconvergent.}\label{table:2}
\end{table} 

\begin{table}[ht] 
\centering
   \begin{tabular}{ c c c c c c } 
  \hline\hline 
   $e_h^q$ & e.c.r. & $ e_h^u$ & e.c.r. &$e_h^{\hat u} $&e.c.r.  \\
    \hline
 2.7400e-002  &  -- & 6.5553e-002   & -- &  1.9182e-002 & -- \\ 
  4.0693e-003  & 2.75 & 1.1591e-002  & 2.50 & 2.4377e-003 & 2.98  \\ 
  5.4030e-004  & 2.91 &  1.6402e-003 & 2.82 & 3.0858e-004 & 2.98 \\
  6.8953e-005  & 2.97 & 2.1698e-004  & 2.92 & 3.8901e-005 & 2.99  
   \\[1ex] \hline \end{tabular} 
\begin{tabular}{ c c c c c c } 
  \hline 
   $\varepsilon_h^u$ & e.c.r. & $ \varepsilon_h^{\hat u}$ & e.c.r. &$\varepsilon_h^\star $&e.c.r.  \\
    \hline
  2.7092e-002  &  -- & 4.7030e-003  & -- & 4.7177e-003 & -- \\ 
  3.1232e-003  & 3.12 & 3.5134e-004 & 3.74 & 3.5580e-004& 3.73  \\ 
  2.1668e-004   &3.85 & 2.2499e-005  & 3.96 & 2.2835e-005 & 3.96 \\
  1.4237e-005  & 3.93 & 1.3859e-006  & 4.02 & 21.4212e-006 & 4.01  
   \\[1ex] \hline \end{tabular} 
\caption{Errors for different triangulations (see Figure \ref{fig:1}) with the lowest order method $k=2$.  All quantities in \eqref{eq:10.2} are shown to be superconvergent. }\label{table:3}
\end{table}

  We finally test the validity of the HDG method as a $p$-method, by fixing the tetrahedrization (the second one in Figure \ref{fig:1}) and increasing $k$ from $0$ to $3$. We compute the relative errors \eqref{eq:10.1} as functions of $k$ and check whether the rates
\begin{equation}\label{eq:10.3}
\frac{\log(e_k/e_{k+1})}{\log(e_{k+1}/e_{k+2})} \approx 1,
\end{equation}   
as would be expected. Note that the theory for $p$-convergence of HDG is not fully developed. The results are reported in Table \ref{table:4}. 
      
\begin{table}[ht] 
\centering
   \begin{tabular}{c c c c c c c } 
  \hline\hline 
 $k$ &  $e_h^q(k)$ & e.c.r. & $ e_h^u(k)$ & e.c.r. &$e_h^{\hat u}(k) $&e.c.r.  \\
    \hline
0 & 4.2490e-001  &  -- & 8.5011e-001  & -- &  3.2710e-001 & -- \\ 
 1 & 3.6794e-002 & -- & 1.2953e-001  & -- & 3.1932e-002 & --  \\ 
2 &  4.0693e-003  & 1.11 &  1.1591e-002 & 0.78 & 2.4377e-003& 0.90 \\
3 &   4.4704e-004  & 1.00 & 1.3590e-003  & 1.13 & 1.7465e-004 & 0.98  
   \\[1ex] \hline \end{tabular} \caption{History of convergence for increasing polynomial degrees, and convergence test following \eqref{eq:10.2}.}\label{table:4} 
\end{table}

Other easy benchmarks for the method are exact polynomial solutions. These have been tried on the implementation as a way to test exactness of approximation and quadrature in the process.

\bibliographystyle{abbrv}
\bibliography{referencesHDGBEM}

\begin{thebibliography}{10}

\bibitem{ArBr:1985}
D.~N. Arnold and F.~Brezzi.
\newblock Mixed and nonconforming finite element methods: implementation,
  postprocessing and error estimates.
\newblock {\em RAIRO Mod\'el. Math. Anal. Num\'er.}, 19(1):7--32, 1985.

\bibitem{BrFo:1991}
F.~Brezzi and M.~Fortin.
\newblock {\em Mixed and hybrid finite element methods}, volume~15 of {\em
  Springer Series in Computational Mathematics}.
\newblock Springer-Verlag, New York, 1991.

\bibitem{ChCo:2012}
Y.~Chen and B.~Cockburn.
\newblock Analysis of variable-degree {HDG} methods for convection-diffusion
  equations. {P}art {I}: general nonconforming meshes.
\newblock {\em IMA J. Numer. Anal.}, 32(4):1267--1293, 2012.

\bibitem{CoGoLa:2009}
B.~Cockburn, J.~Gopalakrishnan, and R.~Lazarov.
\newblock Unified hybridization of discontinuous {G}alerkin, mixed, and
  continuous {G}alerkin methods for second order elliptic problems.
\newblock {\em SIAM J. Numer. Anal.}, 47(2):1319--1365, 2009.

\bibitem{CoGoSa:2010}
B.~Cockburn, J.~Gopalakrishnan, and F.-J. Sayas.
\newblock A projection-based error analysis of {HDG} methods.
\newblock {\em Math. Comp.}, 79(271):1351--1367, 2010.

\bibitem{Ervin:2012}
V.~J. Ervin.
\newblock Computational bases for {$RT_k$} and {$BDM_k$} on triangles.
\newblock {\em Comput. Math. Appl.}, 64(8):2765--2774, 2012.

\bibitem{Felippa:2004}
C.~A. Felippa.
\newblock A compendium of {FEM} integration formulas for symbolic work.
\newblock {\em Engineering Computations}, 21(8):867--890, 2004.

\bibitem{Kirby:2010}
R.~A. Kirby.
\newblock Singularity-free evaluation of collapsed-coordinate orthogonal
  polynomials.
\newblock {\em ACM Trans. Math. Softw.}, 37(1):Article No. 5, 2010.

\bibitem{KiSpCo:2012}
R.~M. Kirby, S.~J. Sherwin, and B.~Cockburn.
\newblock To {CG} or to {HDG}: a comparative study.
\newblock {\em J. Sci. Comput.}, 51(1):183--212, 2012.

\bibitem{Stenberg:1991}
R.~Stenberg.
\newblock Postprocessing schemes for some mixed finite elements.
\newblock {\em RAIRO Mod\'el. Math. Anal. Num\'er.}, 25(1):151--167, 1991.

\bibitem{ZhCuLi:2009}
L.~Zhang, T.~Cui, and H.~Liu.
\newblock A set of symmetric quadrature rules on triangles and tetrahedra.
\newblock {\em J. Comput. Math.}, 27(1):89--96, 2009.

\end{thebibliography}

\end{document}